\documentclass[pdflatex,sn-basic, Numbered]{sn-jnl}
\usepackage{graphicx}%
\usepackage{multirow}%
\usepackage{amsmath,amssymb,amsfonts}%
\usepackage{amsthm}%
\usepackage{mathrsfs}%
\usepackage[title]{appendix}%
\usepackage{xcolor}%
\usepackage{textcomp}%
\usepackage{manyfoot}%
\usepackage{booktabs}%
\usepackage{algorithm}%
\usepackage{algorithmicx}%
\usepackage{algpseudocode}%
\usepackage{listings}%
\usepackage[capitalize, noabbrev]{cleveref}

\theoremstyle{thmstyleone}%
\newtheorem{theorem}{Theorem}[section]
\newtheorem{proposition}[theorem]{Proposition}
\newtheorem{lemma}[theorem]{Lemma}

\theoremstyle{thmstyletwo}%

\newtheorem{problem}[theorem]{Problem}

\theoremstyle{thmstylethree}%
\newtheorem{definition}[theorem]{Definition}

\newcommand{\PG}{\text{PG}}
\newcommand{\AG}{\text{AG}}
\newcommand{\GF}{\text{GF}}

\begin{document}

\title[The arc chromatic number for Galois projective planes, affine planes and Euclidean grids]{The arc chromatic number for Galois projective planes, affine planes and Euclidean grids}

\author[1]{\fnm{Gabriela} \sur{Araujo-Pardo}}\email{garaujo@im.unam.mx}

\author*[2]{\fnm{Leonardo} \sur{Martínez-Sandoval}}\email{leomtz@ciencias.unam.mx}

\affil[1]{\orgdiv{Institute of Mathematics}, \orgname{UNAM}, \city{Juriquilla}, \state{Querétaro}, \country{México}}

\affil*[2]{\orgdiv{Faculty of Sciences}, \orgname{UNAM}, \city{Mexico City}, \state{Mexico City}, \country{Mexico}}

    \abstract{In 1986, Csima and Füredi determined the minimum number of arcs required to partition the points of Galois projective planes $PG(2,q)$ and affine planes $AG(2,q)$. We revisit these foundational results and extend them in several directions.
    
    First, we present alternative proofs for some of their results, providing a more streamlined approach. By leveraging this new perspective, we determine the exact value for a fractional version of the problem in both geometric settings. We also establish a computational framework that yields new constructions and explores the structural properties of the space of colorings.
    
    Finally, we apply these findings to Euclidean grids, improving upon a 2004 construction by Wood. We prove that a partition into $(1+\epsilon)n$ sets in general position exists for any $\epsilon > 0$ and sufficiently large $n$. Additionally, we provide exact minimal partitions for small Euclidean grids.
}


\maketitle

\section{Introduction}

In this work, we study finite planes, which we define as pairs $(\mathcal{P},\mathcal{L})$, where $\mathcal{P}$ and $\mathcal{L}$ denote the sets of points and lines, respectively. Our aim is to understand structural properties of collections of arcs, with particular interest in coloring problems concerning arcs in three specific families: Galois projective planes, finite affine planes, and Euclidean grids.

\subsection{Definitions and Notation}

For each $q > 1$ that is a power of a prime, we denote by \textit{the Galois projective plane $\PG(2,q)$} the finite projective plane over the field $\GF(q)$. The \textit{affine plane $\AG(2,q)$} is obtained from $\PG(2,q)$ by removing one of its lines and all the points on this line. For any positive integer $n$, we define the \textit{Euclidean grid $G_n$} as the finite plane whose point set is $[n] \times [n]$, where $[n] := \{1, \ldots, n\}$, and collinearity is inherited from the Euclidean plane.

We define an \textit{arc} in a finite plane $(\mathcal{P},\mathcal{L})$ as a set $\mathcal{A} \subset \mathcal{P}$ containing no three collinear points\footnote{We focus here on finite planes, which are two-dimensional finite geometries. For $n$-dimensional settings with $n \geq 3$, a set with no three collinear points is called a \textit{cap}, while an \textit{arc} must satisfy the stronger general position condition that no $n+1$ points lie on a common hyperplane.}; that is, for any line $\ell \in \mathcal{L}$ we have $|\ell\cap \mathcal{A}|<3$. We are interested in minimal partitions of the point set of a finite plane into arcs. We pose this as a coloring problem.

For each positive integer $k$, a \textit{$k$-coloring} of $(\mathcal{P},\mathcal{L})$ is a map $c: \mathcal{P} \to [k]$. For a given $k$-coloring and each $j \in [k]$, we refer to $c^{-1}(j)$ as the $j$-th \textit{color class}. We say that a $k$-coloring is \textit{arc-proper} if every color class is an arc.

\begin{definition}
Let $(\mathcal{P}, \mathcal{L})$ be a finite geometry. The \textit{arc chromatic number} $\chi_\mathcal{A}((\mathcal{P},\mathcal{L}))$ is the least $k$ such that there exists an arc-proper $k$-coloring of $(\mathcal{P},\mathcal{L})$.
\end{definition}

We say that two colorings of $\PG(2,q)$ (or $\AG(2,q)$) are \textit{structurally equivalent} if one can be obtained from the other by permuting colors or applying a projective transformation (respectively, an affine transformation). Once we know that an arc-proper $k$-coloring exists for a given finite plane, we can ask whether it is unique modulo structural equivalence. 

Furthermore, we are also interested in a fractional variant of this parameter. A \textit{$(k:b)$-coloring} assigns to each point in $\mathcal{P}$ a set of $b$ colors from $[k]$. In this context, the \textit{$j$-th color class} is the collection of points from $\mathcal{P}$ that received $j$ as one of its colors. We say the $(k:b)$-coloring  is arc-proper when each color class is an arc.

\begin{definition} 
The \textit{$b$-fold arc chromatic number} $\chi_{\mathcal{A},b}((\mathcal{P},\mathcal{L}))$ is the least $k$ such that an arc-proper $(k:b)$-coloring of $(\mathcal{P},\mathcal{L})$ exists. The \textit{fractional arc chromatic number} is defined as: 
$$\chi_{\mathcal{A},f}((\mathcal{P},\mathcal{L}))=\inf_b \frac{\chi_{\mathcal{A},b}((\mathcal{P},\mathcal{L}))}{b}.$$
\end{definition}

These definitions raise the natural question of studying the arc chromatic number and its fractional variant for special families of finite planes.

\subsection{State of the Art}

Coloring problems in projective planes and finite geometries have been extensively studied through many different lenses. A prominent line of research focuses on \textit{rainbow-free colorings}, where the goal is to avoid lines that contain only points of distinct colors. This leads to parameters such as the \textit{upper chromatic number} and the \textit{heterochromatic number} \cite{araujo2003daisy, bacso2013twobloc, bacso2007upper}. Other variations include \textit{balanced} versions of these problems \cite{akm2015balanupper, bbmns21balaupper} and studies on \textit{complete colorings}, which involve the achromatic and pseudoachromatic indices in affine and projective settings \cite{araujo2019indaffine, araujo2019inproy}.

Within this broad landscape, arcs represent a particularly significant class of objects. The motivation for focusing on arcs is rooted in their fundamental structural role in finite geometry. In $\PG(2,q)$, the maximum size of an arc is $q+1$ for odd $q$ and $q+2$ for even $q$. A cornerstone result by Segre \cite{segre1955ovals} established that if $q$ is odd, every oval (an arc of size $q+1$) is a conic. Beyond the plane, determining the size of arcs is a deep problem connected to \textit{maximum distance separable} (MDS) codes and the \textit{MDS conjecture}. For foundational definitions and results on finite projective planes, we refer to the classic text by Hirschfeld \cite{hirschfeld1998projective} and the modern treatment by Kiss and Sz\H{o}nyi \cite{kiss2019finite}. For recent developments on arcs and their connections to codes and the MDS conjecture, we refer to the survey by Ball and Lavrauw \cite{ball2020arcs} and the survey on arcs and generalizations by Hirschfeld and Thas \cite{hirschfeld2025arcs}.

In the specific context of partitioning the point set into arcs, the foundational work of Csima \cite{csima1984colouring} and Csima and Füredi \cite{csima1986colouring} provides the exact values for the arc chromatic number of classical planes.

\begin{theorem}
\label{thm:csimafuredi}
    For any $q$ that is a power of a prime, we have $\chi_\mathcal{A}(\PG(2,q)) = q + 1$. When $q$ is odd, $\chi_\mathcal{A}(\AG(2,q))=q$, and when $q$ is even, $\chi_\mathcal{A}(\AG(2,q))=q-1$.
\end{theorem}

In higher dimensions, the analogous problem for partitioning into caps has also been studied. For example, the study of such partitions for $\PG(d,2)$ is deeply motivated by its connection to Steiner triple systems. In this direction, Fugère, Haddad, and Wehlau \cite{FHW94} proved that $\chi_\mathcal{A}(\PG(5,2)) = 5$, while Haddad \cite{Had99} showed that $\chi_\mathcal{A}(\PG(d,2))$ is unbounded. Other related research considers partitions into even stricter configurations, such as partitions into caps of the same size \cite{ebert1985partitioning} or partitions into normal rational curves \cite{lavrauw2014geometry}.

The problem extends to Euclidean settings through the \textit{no-three-in-line problem}, which asks for the largest arc in a grid $G_n$. For a thorough account of this classic problem, we refer to Chapter F4 in the book by Guy \cite{guy2004unsolved}. Representative contributions to the maximum size of arcs in $G_n$ include work by Hall et al. \cite{hall1975some}, Flammenkamp \cite{flammenkamp1992progress}, and the very recent solution for $n=60$ by Prellberg \cite{prellberg2026constraintsatisfactionprogrammingnothreeinline} using constraint satisfaction programming. Other related results are also discussed in the survey by Chandran, Klavžar and Tuite \cite{v2026generalpositionproblemsurvey}, which covers various general position problems beyond the standard case.

Regarding partitions of the grid, the first bounds for the arc chromatic number were established by Wood \cite{wood2004note}.

\begin{theorem}[Wood, 2004]
    For the Euclidean grid $G_n$, the arc chromatic number satisfies $\chi_{\mathcal{A}}(G_n)\leq 3n-2$. Furthermore, for any $\epsilon > 0$, the bound $\chi_{\mathcal{A}}(G_n)\leq (2+\epsilon)n$ holds asymptotically.
\end{theorem}

Beyond these integer parameters, fractional variants of coloring problems have garnered significant interest within combinatorial optimization and graph theory. Fractional colorings often provide tighter bounds and reveal deeper structural properties of the underlying objects \cite{godsil2001algebraic, scheinerman2011fractional}. However, as far as we are aware, the fractional arc chromatic number of the specific finite geometries investigated in this work has not been previously addressed in the literature.

\subsection{Our Contributions}

In \cref{sec:chromproj}, we revisit the theorem of Csima and Füredi, presenting alternative, more streamlined proofs of some of its components. In particular, we find an optimal arc-proper coloring for $\PG(2,q)$ (and one for $\AG(2,q)$ when $q$ is odd) that is not structurally equivalent to theirs.

\begin{theorem}
\label{thm:nounique} The following holds:
\begin{itemize}
    \item The arc-proper $(q+1)$-coloring of $\PG(2,q)$ is not unique modulo structural equivalence.
    \item The arc-proper $q$-coloring of $\AG(2,q)$, when $q$ is odd, is not unique modulo structural equivalence.
\end{itemize}
\end{theorem}

Building upon this, in \cref{sec:fractional} we determine the exact value of the fractional arc chromatic number for $\PG(2,q)$ and for $\AG(2,q)$.

\begin{theorem}
\label{thm:fractional}
    For any $q$ that is the power of a prime we have
    $$
    \chi_{\mathcal{A},f}(\PG(2,q)) =
    \begin{cases}
        \frac{q^2+q+1}{q+1} & \text{if $q$ is odd}, \\
        \frac{q^2+q+1}{q+2} & \text{if $q$ is even,}
    \end{cases}
    $$

    and 

    $$
    \chi_{\mathcal{A},f}(\AG(2,q)) =
    \begin{cases}
        \frac{q^2}{q+1} & \text{if $q$ is odd}, \\
        \frac{q^2}{q+2} & \text{if $q$ is even.}
    \end{cases}
    $$

\end{theorem}

This result demonstrates a sharp distinction between the regular and fractional chromatic numbers.

In \cref{sec:computational}, we tackle the question of the uniqueness of minimum arc-proper colorings in affine planes of even order. To do so, we present two computational frameworks. We prove the following.

\begin{theorem}
\label{thm:allcomp} The following holds:
\begin{itemize}
    \item There is a unique arc-proper $1$-coloring of $\AG(2,2)$ modulo structural equivalence.
    \item There is a unique arc-proper $3$-coloring of $\AG(2,4)$ modulo structural equivalence.
    \item The arc-proper $7$-coloring of $\AG(2,8)$ is not unique modulo structural equivalence.
\end{itemize}
\end{theorem}

Finally, in \cref{sec:euclidean}, we sharpen the initial bounds established by Wood for the Euclidean grid $G_n$ and provide exact values for small $n$.

\begin{theorem}
    Let $n$ be any positive integer. We have: 
    
    \begin{itemize}
        \item $\frac{n}{2}\leq \chi_{\mathcal{A}}(G_n) \leq 2n$.
        \item For any $\epsilon > 0$, there exists $n_0$ for which  $\chi_{\mathcal{A}}(G_n) \leq (1+\epsilon) n$, when $n\geq n_0$.
        \item The precise values for $\chi_{\mathcal{A}}(G_n)$ for $n \leq 8$ are as follows:
        \begin{center}
        \begin{tabular}{|c|c|c|c|c|c|c|c|c|}
        \hline
        $n$ & 1 & 2 & 3 & 4 & 5 & 6 & 7 & 8 \\
        \hline
        $\chi_{\mathcal{A}}(G_n)$ & 1 & 1 & 2 & 2 & 3 & 4 & 4 & 4\\
        \hline
        \end{tabular}
        \end{center}
    \end{itemize}
\end{theorem}

We conclude in the final section with a discussion and a compilation of open problems.

\section{Revisiting the Standard Arc Chromatic Numbers}
\label{sec:chromproj}

In this section, we discuss \cref{thm:csimafuredi}. We sketch the original proofs given Csima \cite{csima1984colouring} and Csima and Füredi in \cite{csima1986colouring} and provide some alternative proofs that we obtained by studying the problem independently.

When $q$ is odd, a pidgeonhole argument for the following lower bounds dates back to \cite[Theorem 1 and Theorem 2]{csima1984colouring}:
\begin{align*}
  &\chi_\mathcal{A}(\PG(2,q))\geq q+1,\\
  &\chi_\mathcal{A}(\AG(2,q))\geq q.
\end{align*}

For example, for $q$ odd, there are $q^2+q+1$ points in $\PG(2,q)$. An arc-proper coloring using $q$ colors or fewer must yield a color class of size $\left \lceil \frac{q^2+q+1}{q} \right \rceil = q+2$, contradicting the fact that in $\PG(2,q)$ arcs have size at most $q+1$. The second inequality above is proven analogously. In \cite[Theorem 6.1]{csima1986colouring} it was noted that an analogous argument shows that for $q$ even
\begin{align*}
  &\chi_\mathcal{A}(\AG(2,q))\geq q-1.
\end{align*}

For $q$ even, arcs in $\PG(2,q)$ can have size at most $q+2$, so the pidgeonhole argument only yields the weaker bound $\chi_\mathcal{A}(\PG(2,q))\geq q$. In order to improve this value to $q+1$, as needed for \cref{thm:csimafuredi}, a convexity argument is used in \cite[Lemma 2.1] {csima1986colouring}. Their idea is using Jensen's inequality on the sum $\sum_{i=1}^q \binom{\alpha_i}{2}$, where $\alpha_1 \geq \alpha_2 \geq \ldots \geq \alpha_q$ are the sizes of the color classes of a supposed arc-proper $q$-coloring. A double-counting argument then yields a contradiction. When we rediscovered this result, we used a similar approach, but via Karamata's inequality. We present it here for completeness.

\begin{proposition}
    For $q$ a power of $2$, there does not exist an arc-proper $q$-coloring of $\PG(2,q)$.
\end{proposition}

\begin{proof}

Suppose an arc-proper $q$-coloring exists, with chromatic classes $C_1, C_2, \ldots, C_q$ of sizes $\alpha_1 \geq \alpha_2 \geq \ldots \geq \alpha_q$, respectively. We apply Karamata's inequality \cite[Theorem 1]{kadelburg2005inequalities}, paraphrased below for convenience:

\begin{theorem}
    Let $a=(a_1,\ldots,a_n)$ and $b=(b_1,\ldots,b_n)$ vectors in $\mathbb{R}^n$, and $r,s$ real numbers. If the vector $a$ majorizes $b$, i.e.,
    \begin{itemize}
        \item $a_1\geq \ldots \geq a_n$ and $b_1\geq \ldots \geq b_n$,
        \item $a_1+\ldots+a_k \geq b_1+\ldots+b_k$ for each $1\leq k \leq n-1$, and
        \item $a_1+\ldots+a_n=b_1+\ldots+b_n$;
    \end{itemize}
    and $f:(r,s)\to \mathbb{R}$ is a convex function, then $$\sum_{i=1}^n f(a_i)\geq \sum_{i=1}^n f(b_i).$$
\end{theorem}

To use the inequality, we claim that the vector $(\alpha_1, \alpha_2, \ldots, \alpha_q)$ majorizes the vector $(q+2, q+1, \ldots, q+1)$, where $q+1$ appears $q-1$ times. First, note that
$$
\sum_{i=1}^q \alpha_i = \sum_{i=1}^q |C_i| = |\mathcal{P}| = q^2 + q + 1 = q+2 + (q+1)(q-1),
$$
so both vectors have the same sum of entries.

Next, observe that
$$
q\alpha_1 \geq \sum_{i=1}^q \alpha_i = q^2 + q + 1,
$$
so $\alpha_1 \geq \left\lceil \frac{q^2 + q + 1}{q} \right\rceil = q+2$. Since each arc has size at most $q+2$, we have $\alpha_1 = q+2$. This implies $\sum_{i=2}^q \alpha_i = (q^2 + q + 1) - (q+2) = (q-1)(q+1)$. The remaining majorization requirements follow, as the rest of the numbers have $q+1$ as their average. By Karamata's inequality, applied to these vectors and the convex function $\binom{x}{2}$, we obtain
$$
\sum_{i=1}^q \binom{\alpha_i}{2} \geq \binom{q+2}{2} + (q-1)\binom{q+1}{2} = \frac{q^3 + q^2 + 2q + 2}{2}.
$$

We now find an alternative expression for $\sum_{i=1}^q \binom{\alpha_i}{2}$. Each chromatic class $C_i$, being an arc, must span $\binom{\alpha_i}{2}$ distinct lines, so the sum counts the number of pairs $(C_i, \ell)$ where $\ell$ is a line spanned by two points of $C_i$. Alternatively, we count these pairs by fixing the line first. Any line $\ell \in \mathcal{L}$ can have at most $2$ points in each $C_i$, and it can have exactly $2$ points in at most $\left\lfloor \frac{q+1}{2} \right\rfloor = \frac{q}{2}$ of the classes $C_i$. Therefore,
$$
\sum_{i=1}^q \binom{\alpha_i}{2} \leq \frac{q}{2} \cdot |\mathcal{L}| = \frac{q^3 + q^2 + q}{2}.
$$

The two inequalities for $\sum_{i=1}^q \binom{\alpha_i}{2}$ are incompatible. We conclude that $\chi_{\mathcal{A}}(\PG(2,q)) = q$ is impossible for $q$ even as well.

\end{proof}

The proof follows almost verbatim the approach in \cite{csima1986colouring}, only differing in the inequality of convexity that is used. Here the use of either Jensen's or Karamata's inequality yields the same result. However, since Karamata's inequality provides a distinction between integer and non-integer parameters, we expect that for other finite planes it could provide sharper values.

For talking about upper bound constructions, it is convenient to provide a brief reminder on how to coordinatize $\PG(2,q)$ and $\AG(2,q)$. We begin with the Galois projective plane $\PG(2,q)$, whose points are represented by homogeneous coordinates $(x:y:z)$, where $x, y, z \in \GF(q)$ are not all zero. By removing an arbitrary line, typically the line at infinity defined by $z=0$, we obtain the affine plane $\AG(2,q)$. For any point where $z \neq 0$, we can normalize the coordinates to the form $(x/z : y/z : 1)$, or simply $(x : y : 1)$. This allows us to identify the points of this affine plane with the set of ordered pairs $(x, y)$ in $\GF(q)^2$, establishing the standard Cartesian representation of the affine plane. In particular, $\AG(2,q)$ has $q^2$ points.

Now we focus on the upper bounds for \Cref{thm:csimafuredi}, which is the content of \cite[Theorem 5.2]{csima1986colouring}. In order to prove 
\begin{align*}
  &\chi_\mathcal{A}(\PG(2,q))\leq q+1
\end{align*}
the authors focus on zero-sets of expressions of the form $f+cg$ where $c \in \GF(q)$ and $f,g$ are quadratic polynomials over $\GF(q)$. They have to make a very careful selection of polynomials $f$ and $g$ in three possible cases: $q$ odd and $-1$ a quadratic residue in $\GF(q)$, $q$ odd and $-1$ not a quadratic residue in $\GF(q)$ and $q$ even.

Here we present a simpler proof for $PG(2,q)$. For our construction,  we employ the \textit{Singer's cyclic model of $\PG(2,q)$} \cite{singer1938theorem}, in which each point is identified with an element of $\mathbb{Z}_{q^2+q+1}$. In this model, each point possesses an \textit{additive inverse}, and for a set of points $\mathcal{P}$, the \textit{additive inverse} $-\mathcal{P}$ is defined as the set of additive inverses of the points in $\mathcal{P}$. A result of Hall \cite{hall1974difference} establishes that, under these definitions, the additive inverse of any line is an arc; see also \cite{faina2002cyclic} for generalizations.

\begin{proposition}
    Let $q$ be a power of a prime. There exists an arc-proper $(q+1)$-coloring of $\PG(2,q)$.
\end{proposition}

\begin{proof}
    Consider the cyclic model of $\PG(2,q)$. Let $\ell_1, \ldots, \ell_{q+1}$ denote the $q+1$ lines through a fixed point $P$. These lines cover $\PG(2,q)$, so their additive inverses $-\ell_1, \ldots, -\ell_{q+1}$ are arcs that also cover $\PG(2,q)$. Define $c(-P) = 1$, and for each $j \in [q+1]$ and each $P' \in -\ell_j \setminus \{-P\}$, define $c(P') = j$. Since each $-\ell_j$ is an arc, this construction yields an arc-proper coloring of $\PG(2,q)$ using $q+1$ colors.
\end{proof}

Note that this construction is not structurally equivalent to the one given by Csima and Füredi since theirs has $q$ arcs of size $q+1$, and a singleton, while the one above has $q$ arcs of size $q$, and one of size $q+1$.

We will build on top of this construction later on, when we study the fractional variants of the arc chromatic numbers for classical planes.

For $q$ odd, the construction in \cite{csima1986colouring} that shows
\begin{align*}
  &\chi_\mathcal{A}(\AG(2,q))\leq q
\end{align*}
is quite simple: we can partition $\AG(2,q)$ into a set of $q$ conics, namely, using the following parabolas: $\mathcal{P}_i = \{(x, y) : y = x^2 + i\}$ for $i \in \GF(q)$, where we use the coordinatization above. Note that each of these conics has $q$ affine points. The following proposition shows an arc-proper $q$-coloring into one color class of size $1$, and $q-1$ of size $q+1$, and hence a non-structurally equivalent construction.

\begin{proposition}
    Let $q$ be a power of an odd prime. We have
    $$
    \chi_{\mathcal{A}}(\AG(2,q)) \leq q.
    $$
\end{proposition}
 
\begin{proof}
Let $d$ be a quadratic non-residue in $\GF(q)$. For any $c$ in $\GF(q)\setminus \{0\}$, consider $C_c$ as the conic given by the equation
$$f_c(x,y,z):= x^2 - dy^2 - cz^2=0.$$

Because $d$ is a quadratic non-residue, no point with $z=0$ (at infinity) can be on $C_c$, so we may assume $z=1$ above. Since $c\neq 0$, the point $(0:0:1)$ is not in $C_c$. Also, this is an irreducible conic, so no three of its points are collinear. Thus, each $C_c$ is an arc of size $q+1$. Therefore, these $q-1$ conics and the point $(0,0)$ form the desired arc-proper $q$-coloring.
\end{proof}

The construction in \cite[Theorem 6.3]{csima1986colouring} that shows that for $q$ even 
\begin{align*}
  &\chi_\mathcal{A}(\AG(2,q))\leq q-1
\end{align*}
is very similar to the even case for $\PG(2,q)$, in the sense that it relies on finding an appropriate pair of homogeneous quadratic polynomials.

When we independently studied this problem, we only managed to provide an arc-proper $(q-1)$-coloring for $\AG(2,q)$ in the even case for small values of $q$, and we posed as an open problem to find such a construction in general. This led us to consider the computational approaches described in \cref{sec:computational} and it also led to two alternative solutions by Francesco Pavese\footnote{Personal communication} and Péter Sziklai\footnote{Personal communication}. We paraphrase them here as they provide additional points of view on how to construct colorings.

Remember that for $q$ even, we have the following classical result of Qvist \cite{qvist1952some}.

\begin{theorem}
\label{thm:qvist}
Let $q$ be a power of $2$. Let $\mathcal{A}$ be an arc of size $q+1$ in a finite projective geometry of order $q$. Then, there is a unique tangent line through each point of $\mathcal{A}$, and all of these tangents pass through a common point $P$ outside $\mathcal{A}$.
\end{theorem}

We call the point $P$ in Theorem \ref{thm:qvist} the \textit{nucleus} of $\mathcal{A}$. The set $\mathcal{A} \cup \{P\}$ is still an arc. When we have a conic in $\PG(2,q)$ given by a quadratic equation, the nucleus can be found as the point that vanishes the partial derivatives.

\begin{proposition}
\label{prop:dobleprueba}
    Let $q$ be a power of $2$. There exists an arc-proper $(q-1)$-coloring of $\AG(2,q)$.
\end{proposition}

\begin{proof} (by Francesco Pavese)
Fix $e$ in $\GF(q)$ such that $x^2+x+e$ is irreducible over $\GF(q)$. We will set $z=0$ as the line at infinity. For each non-zero $a$ in $\GF(q)$, consider $C_a$ the conic with equation $$f_a(x,y,z):=x^2+xy+ey^2+az^2=0.$$

None of these conics have a point at infinity, since $x^2+xy+ey^2$ has no solution over $\GF(q)$. Also, since $a\neq 0$, then $(0:0:1)$ is not on any of these conics. Moreover, these conics are disjoint, since a common solution $(x:y:z)$ in $C_a$ and $C_a'$ implies
\begin{align*}
    x^2+xy+ey^2+az^2&=0\\
    x^2+xy+ey^2+a'z^2&=0,
\end{align*}
and by subtracting them we get $(a-a')z^2=0$, which implies $a=a'$.

Therefore, these $q-1$ conics are a partition for the $q^2-1=(q-1)(q+1)$ points of $\AG(2,q)\setminus \{(0:0:1)\}$. Now, we note that the partial derivatives of $f_a$ are (using that $GF(q)$ has characteristic $2$):
$$\frac{\partial f}{\partial x} = y\quad \frac{\partial f}{\partial y} = x \quad \frac{\partial f}{\partial z} = 0.$$

For any $a$, they vanish exactly when $x=y=0$, that is, in the point $(0:0:1)$. Therefore all of these conics are non-degenerate and $(0:0:1)$ is precisely the nucleus of each of them. By adding this point to any of the conics, we obtain the desired coloring.

\end{proof}

\begin{proof} (by Péter Sziklai) The affine plane $\AG(2,q)$ and $\GF(q^2)$ can be identified as $2$-dimensional vector spaces over $\GF(q)$. The origin will be $0\in \GF(q^2)$, and the lines through it are the sets \[ L_a=\{ax:x\in \GF(q)\}, \] for any $a\in \GF(q^2)^*$. The points $S=\{s\in \GF(q^2):s^{q+1}=1\}$ form a conic (the ``unit circle''), whose intersection with a line is \[ S\cap L_a = \{ax:(ax)^{q+1}=1,\ x\in \GF(q)\}. \] But \[ 1=(ax)^{q+1} =a^{q+1}x^{q+1} =a^{q+1}x^2, \] which has a unique solution $x\in \GF(q)$ when $q$ is even. Hence every line through $0$ meets $S$ in exactly one point, so $0$ is the nucleus of $S$. It also follows that, by magnifying $S$ from the origin by $b\in \GF(q)^*$, we obtain $q-1$ pairwise disjoint conics $bS=\{bs:s\in S\},$ each having nucleus $0$. Altogether, they partition the $(q-1)(q+1)$ points of $\AG(2,q)\setminus\{0\}$. Finally, assigning color $b$ to the conic $bS$, and adding the origin to any one of the color classes, we obtain an arc-proper $(q-1)$-coloring of $\AG(2,q)$. \end{proof}

Under the identification $\GF(q^2)=\GF(q)(\alpha)$, where
$\alpha^2+\alpha+e=0$, the norm of $x+y\alpha$ is $x^2+xy+ey^2$. Consequently, the sets $bS$ in the second proof are precisely the affine
conics $x^2+xy+ey^2=b^2$. 
Since squaring is a bijection of $\GF(q)^*$ when $q$ is even, these are
exactly the conics $x^2+xy+ey^2=a$ appearing in the first proof.

Another perspective on the coloring  of $\AG(2,4)$ exploiting the so-called \textit{Daisy Structure} is presented in Appendix \ref{app:daisy}. Also, the computational approach in \cref{sec:computational}  provides structurally distinct colorings for $\AG(2,8)$ that complement this viewpoint.

\section{The Fractional Arc Chromatic Numbers}
\label{sec:fractional}

We now focus on proving \cref{thm:fractional}. For a finite plane $(\mathcal{P},\mathcal{L})$, let $\alpha((\mathcal{P},\mathcal{L}))$ denote the size of a maximal arc. The lower bounds below all follow from the inequality 
\begin{align*}
    \chi_{\mathcal{A},f}((\mathcal{P},\mathcal{L})) \geq \frac{|\mathcal{P}|}{\alpha((\mathcal{P},\mathcal{L}))}.
\end{align*}

The proof of this inequality is standard. We provide it for completeness. Since $\chi_{\mathcal{A},f}$ is an infimum, it suffices to prove for each $b$ that
\begin{align*}
    \frac{\chi_{\mathcal{A},b}((\mathcal{P},\mathcal{L}))}{b} \geq \frac{|\mathcal{P}|}{\alpha((\mathcal{P},\mathcal{L}))},
\end{align*}
i.e., $\chi_{\mathcal{A},b}((\mathcal{P},\mathcal{L})) \cdot \alpha((\mathcal{P},\mathcal{L})) \geq b \cdot |\mathcal{P}|$. This is clear since the sum of the cardinalities of the color classes is bounded by $\chi_{\mathcal{A},b}\cdot\alpha((\mathcal{P},\mathcal{L}))$, and this is at least $b \cdot |\mathcal{P}|$ because the color classes cover each point $b$ times.

Then, what \cref{thm:fractional} claims is that this lower bound is always attained. We provide all the necessary constructions below. Note that a transitivity argument on the symmetry group would also yield these constructions. We prefer an explicit approach here, particularly to tie the constructions to specific curves and because Lemma \ref{lem:bijection} below is of independent interest.

\begin{proposition}
    Let $q$ be a power of an odd prime. We have
    $$
    \chi_{\mathcal{A},f}(\PG(2,q)) \leq \frac{q^2+q+1}{q+1}.
    $$
\end{proposition}

\begin{proof}
    We again use Singer's cyclic model of $\PG(2,q)$. As color classes, we take the additive inverse of every line, which are known to be arcs. There are $q^2+q+1$ such color classes. Each point lies on $q+1$ lines, so this yields an arc-proper $(q^2+q+1 : q+1)$-coloring. Therefore, $\frac{q^2+q+1}{q+1}$ is one of the values considered in the definition of $\chi_{\mathcal{A},f}$ as an infimum, and the claimed inequality holds.
\end{proof}

The construction for $q$ even is similar, but we require an auxiliary result that states that, in Singer's cyclic model of $\PG(2,q)$, the arcs obtained as additive inverses of distinct lines have distinct nuclei.

\begin{lemma}
\label{lem:bijection}
    Let $q$ be a power of $2$. Consider the cyclic model of $\PG(2,q)$ and distinct lines $\ell_1$ and $\ell_2$. Then the nuclei of $-\ell_1$ and $-\ell_2$ are distinct. Therefore, every point is the nucleus of exactly one arc of the form $-\ell$.
\end{lemma}

\begin{proof}
    Suppose, for contradiction, that the nuclei of both $-\ell_1$ and $-\ell_2$ are the same point $Q$. The lines $\ell_1$ and $\ell_2$ intersect in a single point $P$. Let $\ell_3, \ldots, \ell_{q+1}$ be the remaining lines through $P$. Then, the arcs $-\ell_1, \ldots, -\ell_{q+1}$ intersect pairwise in the single point $-P$. By Theorem \ref{thm:qvist}, each of these $q+1$ arcs has exactly one tangent through $-P$. Since the line through $Q$ and $-P$ is the tangent to both $-\ell_1$ and $-\ell_2$, one of the $q+1$ lines through $-P$, say $m$, is not the tangent of any of the arcs $-\ell_1, \ldots, -\ell_{q+1}$. But then $m$ must contain the point $-P$ and at least one other point in each of these arcs, yielding the contradiction $|m| \geq q+2$.

    This shows that the map assigning to each line the nucleus of its additive inverse is injective. Since the number of lines and points is both $q^2+q+1$, this map is bijective.
\end{proof}

We are now ready to provide a coloring for the case when $q$ is even.

\begin{proposition}
    Let $q$ be a power of $2$. We have
    $$
    \chi_{\mathcal{A},f}(\PG(2,q)) \leq \frac{q^2+q+1}{q+2}.
    $$
\end{proposition}

\begin{proof}
    Consider the cyclic model of $\PG(2,q)$. For each line $\ell$, let $P_\ell$ be the nucleus of $-\ell$. We propose as color classes the sets $-\ell \cup \{P_\ell\}$ as $\ell$ ranges over the lines of $\PG(2,q)$. These color classes are arcs of size $q+2$. Each point lies on $q+1$ lines, so it is in $q+1$ arcs $-\ell$; and by Lemma \ref{lem:bijection}, it is the nucleus of an arc $-\ell$ exactly once. Therefore, each point is in exactly $q+2$ classes.

    Thus, we have constructed an arc-proper $(q^2+q+1 : q+2)$-coloring, from which
    $$
    \chi_{\mathcal{A},f} \leq \frac{q^2+q+1}{q+2},
    $$
    as desired.
\end{proof}

To provide fractional colorings for the affine plane, we use conics with specific equations.

\begin{proposition}
    Let $q$ be a power of an odd prime. We have
    $$
    \chi_{\mathcal{A},f}(\AG(2,q)) \leq \frac{q^2}{q+1}.
    $$
\end{proposition}
 
\begin{proof}
Let $d$ be a quadratic non-residue in $\GF(q)$. For any $a,b$ in $\GF(q)$, consider $C_{a,b}$ as the conic given by the equation
$$f_{a,b}(x,y,z):= (x-az)^2 - d(y-bz)^2 - z^2=0$$

Because $d$ is a quadratic non-residue, no point with $z=0$ (at infinity) can be on this conic. Also, this is an irreducible conic, so no three of its points are collinear. Thus, each $C_{a,b}$ is an arc of size $q+1$. As there are $q$ choices for $a$ and $q$ choices for $b$, the family $\mathcal{F}$ contains exactly $q^2$ such arcs.

To determine the number of arcs passing through a fixed affine point $(x_0, y_0) \in AG(2,q)$, we count the pairs $(a,b)$ satisfying:
$$(x_0 - a)^2 - d(y_0 - b)^2 = 1.$$

For $x_0$ and $y_0$ fixed, this is now a non-degenerate conic on variables $a,b$, so it has precisely $q+1$ solutions. Therefore, every point in $AG(2,q)$ is covered exactly $q+1$ times by the conics $C_{a,b}$. This yields an arc-proper $(q^2:q+1)$-coloring, as desired.
\end{proof}

\begin{proposition}
    Let $q$ be a power of two. We have
    $$
    \chi_{\mathcal{A},f}(\AG(2,q)) \leq \frac{q^2}{q+2}.
    $$
\end{proposition}

\begin{proof}
    Let $e \in \GF(q)$ be such that $x^2+x+e$ is irreducible. For any $a, b \in \GF(q)$, consider the conic $C_{a,b}$ defined by the equation:
    $$f_{a,b}(x,y,z)= ex^2 + xy + y^2 + bxz + ayz + (1 + ea^2 + ab + b^2) z^2 = 0.$$

    Setting $z=0$ yields $ex^2 + xy + y^2 = 0$, which has no non-zero solutions since $x^2+x+e$ is irreducible. Thus, $C_{a,b}$ has no points at infinity. In characteristic $2$, the nucleus $N_{a,b}$ is found by setting the formal partial derivatives to zero:
    $$\frac{\partial f}{\partial x} = y + bz = 0, \quad \frac{\partial f}{\partial y} = x + az = 0.$$
    Dehomogenizing with $z=1$ yields the nucleus $N_{a,b} = (a,b)$. Evaluating $f_{a,b}$ at $(a,b,1)$ gives $ea^2 + ab + b^2 + ab + ab + (1 + ea^2 + ab + b^2) = 1 \neq 0$, so the conic is non-degenerate. 

    Defining $H_{a,b} = C_{a,b} \cup \{(a,b)\}$ gives a family of $q^2$ hyperovals in $AG(2,q)$, each being an arc of size $q+2$. A fixed point $(x_0,y_0)$ is the nucleus of $H_{a,b}$ if $(a,b) = (x_0, y_0)$ (one case). It lies on the conic $C_{a,b}$ if:
    $$ea^2 + ab + b^2 + y_0 a + x_0 b + (ex_0^2 + x_0 y_0 + y_0^2 + 1) = 0.$$
    This is an irreducible conic on variables $a,b$ with no points at infinity, so it has exactly $q+1$ solutions for $(a,b)$. In total, each point $(x_0, y_0)$ is covered by $1 + (q+1) = q+2$ arcs, providing the required arc-proper $(q^2 : q+2)$-coloring.
\end{proof}

\section{Computational Frameworks for Arc-Proper Colorings}
\label{sec:computational}

Here we present two distinct computational frameworks to tackle arc-coloring problems. Within each of them, we first developed a general-purpose finite field (GF) factory for arithmetic in $\GF(2^k)$ for $k \in \{2, 3, 4\}$. This factory employs explicit irreducible polynomials and supports efficient multiplication via binary polynomial arithmetic. For an introduction to the computational aspects of this approach, see Part 5 of \cite{arndt2010matters}.

Using this arithmetic foundation, we systematically generated the projective points of $\PG(2, 2^k)$ by normalizing homogeneous coordinates $(x:y:z)$ to ensure a canonical representative for each equivalence class, and then extracted the affine points as those with $z\neq 0$. Projective lines were constructed by enumerating normalized coefficient triples and aggregating points satisfying the corresponding linear equations. A determinant-based collinearity test using $\GF(2^k)$ arithmetic was implemented. 

To verify the correctness of the colorings, we also implemented a slow but easily verifiable $O(n^3)$ determinant test to determine whether a set of $n$ points in $\PG(2,q)$ forms an arc.

This section is complemented by the experimental scripts available at the GitHub repository \cite{Martinez2026}. We recommend consulting the source code for specific implementation details regarding data structures and search heuristics.

\subsection{Framework 1: Backtracking + forward-checking + MRV}

Our first approach is to try to construct arc-proper colorings by trial and error using a backtracking procedure that tries to incrementally color $n$ points using a fixed amount $k$ of colors. The principal challenge in finding such a coloring with this approach is navigating the immense combinatorial search space. For example, finding an arc-proper $7$-coloring for $\AG(2,8)$ has a $7^{64}$ search space.

To address this, we implemented a backtracking solver enhanced by forward-checking and the Minimum Remaining Values (MRV) heuristic, classical techniques pioneered by Haralick and Elliott \cite{haralick1980increasing}. By performing early look-aheads and prioritizing the most constrained points, the algorithm effectively prunes the search tree and avoids the exploration of non-viable configurations. Despite the potential for exponential growth, the implementation proved remarkably efficient, converging to the coloring in the following proposition in approximately ten seconds on a standard personal workstation.

\begin{proposition}
\label{prop:computational}
    There exists an arc-proper $7$-coloring of $\AG(2,8)$ with color classes of size $10,10,9,9,9,9,8$. Therefore, this coloring is not structurally equivalent to the ones in \cite[Theorem 6.3]{csima1986colouring} or \Cref{prop:dobleprueba}. 
\end{proposition}

\begin{proof} 

  We explicitly present the coloring, beginning with a clarification of our notation. We construct $\GF(8)$ as the set of polynomials over $\GF(2)$, modulo the irreducible polynomial $x^3+x+1$. Each integer $j$ in $\{0,1,\ldots,7\}$ is identified with an element of $\GF(8)$ by interpreting the binary digits of $j$ as coefficients. For example, $5$ in binary is $101$ and thus corresponds to the polynomial $x^2+1$.

  The arc-proper $7$-coloring is given by the following chromatic classes:

  \begin{align*}
      &C_1=\{(0, 0), (0, 1), (1, 0), (1, 1), (7, 5), (3, 5), (3, 7), (5, 7), (5, 3), (7, 3)\}\\
      &C_2=\{(0, 5), (0, 6), (5, 0), (6, 0), (2, 2), (4, 4), (5, 2), (4, 6), (6, 4), (2, 5)\}\\
      &C_3=\{(0, 7), (0, 2), (7, 0), (1, 3), (4, 7), (7, 1), (2, 3), (2, 1), (4, 2)\}\\
      &C_4=\{(0, 3), (0, 4), (4, 0), (6, 6), (7, 7), (6, 7), (4, 3), (1, 6), (7, 4)\}\\
      &C_5=\{(5, 5), (3, 3), (5, 1), (2, 4), (1, 5), (6, 3), (2, 7), (3, 1), (1, 7)\}\\
      &C_6=\{(2, 0), (3, 0), (7, 2), (2, 6), (1, 4), (6, 5), (7, 6), (3, 2)\}\\
      &C_7=\{(1, 2), (3, 6), (5, 4), (6, 1), (3, 4), (6, 2), (4, 5), (5, 6), (4, 1)\}\\
  \end{align*}

  The correctness of this coloring can be verified by a computational determinant test.
\end{proof}

\subsection{Framework 2: Enumeration of arcs + Partition}

Framework 1 is sufficient when we want to find a single coloring with specific properties. However, when we want to study all possible colorings, it becomes redundant, as many large arcs get rediscovered during the backtracking.

Another approach is to first enumerate all possible arcs of interest, and then use a partition (or exact cover) algorithm.

Consider, for example, finding all possible arc-proper $3$-colorings of $\AG(2,4)$. We can proceed as follows:

\begin{enumerate}
    \item Since an arc can have size at most $6$, we can only have colorings with color classes of sizes $6,6,4$ or $6,5,5$.
    \item We enumerate all possible arcs of sizes $4,5,6$.
    \item We use a specialized algorithm to check all the possible choices of three of the enumerated arcs that partition $\AG(2,4)$.
\end{enumerate}

For the last part, we typically use Knuth's Dancing Links algorithm for exact covers, see \cite{knuth2000dancing}. By executing this exact cover approach on $\AG(2,4)$, we can exhaustively classify its colorings and obtain the following uniqueness result.

\begin{proposition}
\label{prop:unique}
    There exists a unique arc-proper $3$-coloring of $\AG(2,4)$, modulo structural equivalence. 
\end{proposition}

\begin{proof}
    We compute all possible colorings as described above and obtain $288$ distinct valid colorings. To test for equivalence, we consider the group of $2880$ affine transformations, viewed as projective transformations fixing the line $z=0$. Applying these transformations to our set of colorings reveals that the entire set forms a single orbit. Therefore, all these colorings are structurally equivalent, which establishes uniqueness.
\end{proof}

\begin{proof}[Proof of \Cref{thm:allcomp}]
The theorem claims results on the uniqueness of minimum colorings for $AG(2,q)$. The uniqueness in the case $q=2$ is immediate. The case $q=4$ follows by \Cref{prop:unique}. The existence of non-structurally equivalent colorings in the case $q=8$ follows from \Cref{prop:computational} and either \cite[Theorem 6.3]{csima1986colouring} or \Cref{prop:dobleprueba}.
\end{proof}

\section{The Arc Chromatic Number for Euclidean Grids}
\label{sec:euclidean}

In the Euclidean grid $G_n$, an arc can contain at most two points per row, and thus at most $2n$ points in total. Consequently, a pigeonhole argument implies that the arc chromatic number satisfies $\chi_{\mathcal{A}}(G_n) \geq \left \lceil \frac{n^2}{2n}\right \rceil = \left \lceil \frac{n}{2} \right \rceil$. We now establish a general upper bound for the arc chromatic number of Euclidean grids.

\begin{proposition}
\label{prop:primes}
    For any positive integer $n$, the grid $G_n$ admits an arc-proper $2n$-coloring.
\end{proposition}

\begin{proof}
    By Bertrand's postulate, there exists a prime number $q$ such that $n \leq q \leq 2n$. We identify $G_n$ as a sub-grid of $\AG(2,q)$ and utilize the arc-proper $q$-coloring of $\AG(2,q)$. Since any three collinear points in the Euclidean grid $G_n$ are also collinear in $\AG(2,q)$, the arc-proper coloring of $\AG(2,q)$ induces an arc-proper coloring of $G_n$ using at most $q \leq 2n$ colors.
\end{proof}

Asymptotically, improved versions of Bertrand's postulate yield tighter bounds on the number of colors required. Specifically, for any $\epsilon > 0$, there exists an integer $n_0$ such that for all $n \geq n_0$, there exists a prime $q$ satisfying $n \leq q \leq (1+\epsilon)n$ (see, for instance, \cite{hardy75}). Therefore, for $n \geq n_0$, the arc chromatic number of the Euclidean grid $G_n$ is bounded above by $(1+\epsilon)n$.

We also employed computational methods to determine the arc chromatic number for small Euclidean grids. The results and the techniques used are presented below.

\begin{proposition}
\label{prop:euclides}
    For $n\leq 8$, the precise values for $\chi_{\mathcal{A}}(G_n)$ are given below:

\begin{center}
\begin{tabular}{|c|c|c|c|c|c|c|c|c|}
\hline
$n$ & 1 & 2 & 3 & 4 & 5 & 6 & 7 & 8 \\
\hline
$\chi_{\mathcal{A}}(G_n)$ & 1 & 1 & 2 & 2 & 3 & 4 & 4 & 4\\
\hline
\end{tabular}
\end{center}
\end{proposition}

\begin{proof}
For $n = 1, 2$, all points can be assigned the same color. For $n \geq 3$, we explicitly provide the color classes. Since $G_{n-1}$ is a subgrid of $G_n$, it is enough to provide the colorings for $n=4,5,8$.

$n = 4$
\begin{align*}
  &C_1 = \{(1, 1), (1, 2), (2, 3), (2, 4), (3, 1), (3, 2), (4, 3), (4, 4)\},\\
  &C_2 = \{(1, 3), (1, 4), (2, 1), (2, 2), (3, 3), (3, 4), (4, 1), (4, 2)\}.
\end{align*}

$n = 5$
\begin{align*}
  &C_1 = \{(1, 1), (1, 2), (2, 1), (2, 4), (3, 4), (3, 5), (4, 2), (4, 3), (5, 3), (5, 5)\},\\
  &C_2 = \{(1, 3), (1, 4), (2, 3), (2, 5), (3, 1), (4, 4), (4, 5), (5, 1), (5, 2)\},\\
  &C_3 = \{(1, 5), (2, 2), (3, 2), (3, 3), (4, 1), (5, 4)\}.
\end{align*}

$n = 8$, see also \cref{fig:8x8coloring}.

\begin{figure}
    \centering
    \includegraphics[width=0.4\linewidth]{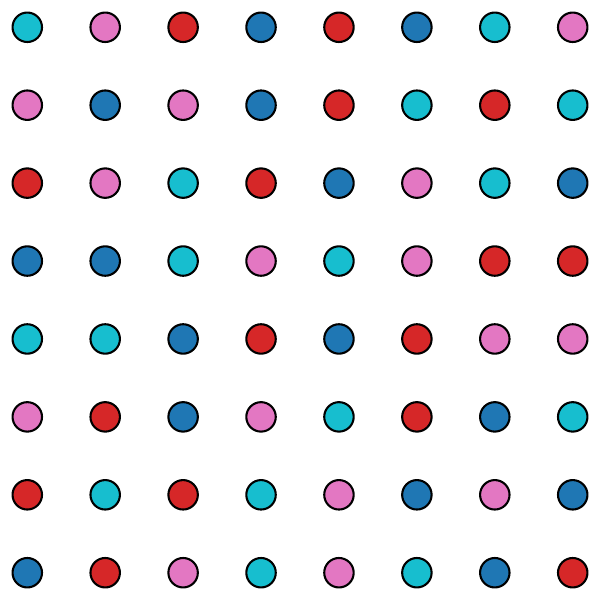}
    \caption{An arc-proper $4$-coloring of $G_8$.}
    \label{fig:8x8coloring}
\end{figure}

\begin{align*}
    C_1 = \{&(7, 3), (1, 5), (8, 2), (3, 3), (1, 1), (5, 4), (4, 8), (6, 2),\\ 
    &(2, 5), (6, 8), (3, 4), (5, 6), (4, 7), (7, 1), (2, 7), (8, 6)\},\\
    C_2 = \{&(1, 2), (8, 5), (2, 3), (3, 2), (3, 8), (5, 7), (7, 7), (7, 5), \\
    &(8, 1), (2, 1), (4, 4), (1, 6), (6, 4), (5, 8), (4, 6), (6, 3)\},\\
    C_3 = \{&(5, 1), (6, 6), (8, 8), (4, 5), (2, 6), (7, 2), (6, 5), (3, 1),\\ 
    &(1, 7), (8, 4), (1, 3), (2, 8), (3, 7), (4, 3), (7, 4), (5, 2)\},\\
    C_4 = \{&(5, 5), (1, 8), (3, 5), (7, 6), (4, 2), (2, 2), (1, 4), (5, 3),\\ 
    &(4, 1), (7, 8), (8, 7), (6, 1), (6, 7), (8, 3), (3, 6), (2, 4)\}.
\end{align*}

The validity of these colorings and the non-existence of arc-proper colorings with fewer colors can be verified computationally using the accompanying code or independently.
\end{proof}

A suitable adaptation of Framework 1 was enough to derive the colorings for $n \leq 7$ and to confirm that no arc-proper colorings with fewer colors exist.

However, for $n = 8$, Framework 1 exceeded our computational limits. To address this, we adapted Framework 2 as follows.

In $G_8$, an arc contains at most $16$ points (at most two points per row). Consequently, in an arc-proper $4$-coloring of $G_8$, each color class must contain exactly $16$ points. For an arc to have $16$ points, it must have two points on each row. This allows a fast recursive enumeration of all possible arcs of size $16$ in $G_8$. After five minutes of computation on a standard personal workstation, we identified all $380$ of them. Subsequently, a backtracking solver was employed to determine whether four of these arcs could form a partition of $G_8$. A valid coloring was found almost instantaneously.

\section{Discussion and Open Problems}

The determination of the arc chromatic numbers for $G_n$ remains an open question, as originally posed by Wood \cite{wood2004note}.

\begin{problem}
Let $n\geq 9$ be a positive integer. Determine the value of $\chi_{\mathcal{A}}(G_n)$.
\end{problem}

Given the inherent complexity of the no-three-in-line problem, we anticipate that this problem will also be challenging. 

After the first version of this work, we were contacted by two independent teams concerning partial progress. First, Haoxuan (Jason) Dong sent us\footnote{Personal communication} a $5$-coloring of $G_9$, which by the lower bound above, would settle  \[ \chi_{\mathcal A}(G_9)=5.\] Dong subsequently joined Thomas Prellberg, Ivet Lobo, and Matthew Lewis to continue investigating this and related grid-coloring problems.

We were also contacted\footnote{Personal communication} by Youngwoo Yoon and Daegeun Lee (Superintelligence Labs Inc., Seoul, Republic of Korea), together with Donghyun Lim and Seungwoo Schin (Mazigle Inc., Seoul, Republic of Korea). They also independently found a $5$-coloring of $G_9$, and moreover communicated that they had established \[ \chi_{\mathcal A}(G_{10})=6. \]

In both cases, we were able to independently verify that the colorings are arc-proper. Both teams are preparing manuscripts with further details and verification materials.

The fractional arc chromatic number can also be defined for $G_n$.

\begin{problem}
Let $n$ be a positive integer. Determine $\chi_{\mathcal{A},f}(G_n)$.
\end{problem}

The approach from \cref{prop:primes} shows that $$\chi_{\mathcal{A},f}(G_n)\leq \frac{4n^2}{2n+1},$$
but we have not tried to get a better bound.

\Cref{thm:nounique} and \Cref{thm:allcomp} suggest that the space of arc-proper colorings can be quite large and structurally rich. It raises a natural question.

\begin{problem}
For $q$ a power of a prime, and modulo structural equivalence:
\begin{itemize}
    \item Classify all the arc-proper $\chi_\mathcal{A}(\PG(2,q))$-colorings of $\PG(2,q)$.
    \item Classify all the arc-proper $\chi_\mathcal{A}(\AG(2,q))$-colorings of $\AG(2,q)$.
\end{itemize}
\end{problem}

In particular, our computational framework can not currently determine whether the minimum coloring for $\AG(2,q)$ is structurally unique for $q=2^r$ when $r\geq 4$.

The study of the arc chromatic number can also be extended to other finite planes or finite geometries. We believe this represents a promising direction for future research.
\begin{appendices}
\section{A Daisy Structure Perspective}
\label{app:daisy}

This appendix presents an alternative viewpoint on certain colorings of $\PG(2,q)$ when $q$ is even, commonly referred to as the \textit{Daisy Structure} (Araujo-Pardo \cite{araujo2003daisy}). This perspective may generalize and can be useful when designing or guiding constructions in small even orders such as $q=4$.

\begin{figure}
    \centering
    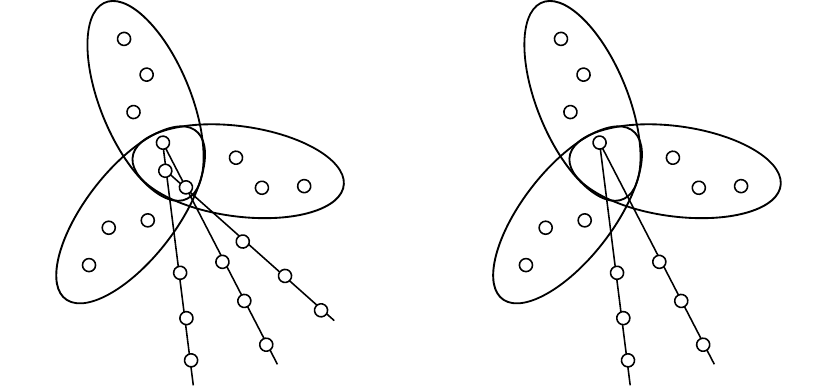
    \caption{The Daisy Structure of $\PG(2,4)$ (left) and the induced structure on $\AG(2,4)$ (right).}
    \label{fig:daisies}
\end{figure}

When $q$ is even, $\PG(2,q)$ admits the following structure. There exists a family of $q-1$ hyperovals (arcs of size $q+2$), denoted $\mathcal{H}_1,\ldots,\mathcal{H}_{q-1}$, whose common intersection is the set $\mathcal{C}=\{(1:0:0),(0:1:0),(0:0:1)\}$, called the \textit{Daisy Center}. Each petal is $\mathcal{P}_i=\mathcal{H}_i\setminus\mathcal{C}$ ($i=1,\ldots,q-1$). The remaining points lie on the three lines through the center, called the \textit{stems} and denoted $\mathcal{S}_0,\mathcal{S}_1,\mathcal{S}_2$. In this framework, lines fall into three types: (1) stems through two center points; (2) lines meeting the center in one point, meeting the opposite stem and each petal in one point; (3) lines meeting no center point, meeting $\frac{q-2}{2}$ petals in two points each and meeting each stem in one non-center point. See \cite{araujo2003daisy} for details; Figure \ref{fig:daisies} illustrates the case $q=4$.

Removing a line (for instance $z=0$) yields a decomposition of $\AG(2,4)$ into a single center point, three petals, and two stems. For $\AG(2,4)$ the petals and stems can be written explicitly as
\begin{align*}
\mathcal{P}_0&=\{(1,1),(2,3),(3,2)\},\\
\mathcal{P}_1&=\{(1,3),(2,2),(3,1)\},\\
\mathcal{P}_2&=\{(1,2),(2,1),(3,3)\},\\
\mathcal{S}'_1&=\{(0,0),(1,0),(2,0),(3,0)\},\\
\mathcal{S}'_2&=\{(0,0),(0,1),(0,2),(0,3)\}.
\end{align*}

\begin{figure}
    \centering
    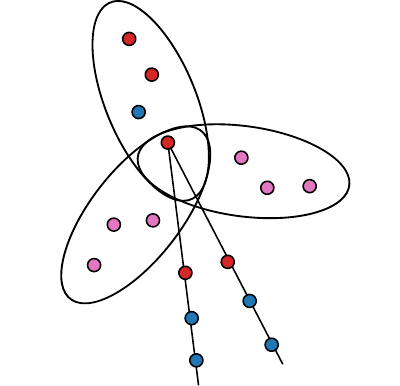
    \caption{An arc-chromatic $3$-coloring for $\AG(2,4)$ arising from the Daisy Structure.}
    \label{fig:daisy-colors}
\end{figure}

An arc-proper $3$-coloring for $\AG(2,4)$ is obtained by grouping petals and stems as follows (see Figure \ref{fig:daisy-colors}):
\begin{align*}
      &C_1=\{(1, 1), (2, 3), (3, 2), (1, 3), (2, 2), (3,1)\},\\
      &C_2=\{(0, 0), (1, 2), (2, 1), (1, 0), (0, 1)\},\\
      &C_3=\{(0, 2), (0, 3), (2, 0), (3, 0), (3, 3)\}.
\end{align*}

Here $C_1=\mathcal{P}_0\cup\mathcal{P}_1$ is the union of two petals; lines of each type intersect these sets in at most two points, and similar checks show $C_2$ and $C_3$ contain no three collinear points. The geometric partition above therefore produces the stated arc-proper coloring.

\end{appendices}

\backmatter

\bmhead{Acknowledgements}

The authors thank Dr. Kolja Knauer for some  suggested bibliography and the discussions shared during his research stay at the Institute of Mathematics in Juriquilla in July 2025.

We also would like to thank several colleagues for the valuable feedback and insights provided following the initial announcement of a previous version of this work. We are grateful to Francesco Pavese for his construction for \Cref{prop:dobleprueba} and for providing additional references on partitions of projective spaces into arcs and caps.  We also thank Péter Sziklai for independently providing another construction for \Cref{prop:dobleprueba} and for various helpful comments that improved the presentation of this work. We thank James Tuite for sharing several relevant references concerning the no-three-in-line problem, in particular, for pointing us toward the work of Wood \cite{wood2004note}, where the coloring problem for the Euclidean variant was first considered. Finally, we thank Haoxuan (Jason) Dong, who  contacted us with results concerning the arc chromatic numbers of the grids $G_n$, and Youngwoo Yoon, Daegeun Lee, Donghyun Lim, and Seungwoo Schin, who subsequently contacted us independently with further results on this problem.

This work was supported by UNAM-PAPIIT IN113324, by UNAM-PAPIIT IN119026 and by SECIHTI-M\'exico CBF2023-2024-552.

\section*{Declarations}

The authors have no relevant financial or non-financial interests to disclose. All authors contributed to the study of the problems in this work. All authors participated in writing the draft of the manuscript, commented on its versions and approved its final version.

\bibliography{refs.bib}

@article{csima1986colouring,
  title={Colouring finite incidence structures},
  author={Csima, J and F{\"u}redi, Zolt{\'a}n},
  journal={Graphs and Combinatorics},
  volume={2},
  number={1},
  pages={339--346},
  year={1986},
  publisher={Springer}
}

@article{segre1955ovals,
  title={Ovals in a finite projective plane},
  author={Segre, Beniamino},
  journal={Canadian Journal of Mathematics},
  volume={7},
  pages={414--416},
  year={1955},
  publisher={Cambridge University Press}
}

@misc{Martinez2026,
  author = {Martínez-Sandoval, Leonardo and Araujo-Pardo, Gabriela},
  title = {Coloring the {E}uclidean Grid and Affine Planes},
  year = {2026},
  publisher = {GitHub},
  journal = {GitHub repository},
  howpublished = {\url{https://github.com/leomtz/coloracion-arc-proyectivo}}
}

@article{ball2020arcs,
  title={Arcs in finite projective spaces},
  author={Ball, Simeon and Lavrauw, Michel},
  journal={EMS Surveys in Mathematical Sciences},
  volume={6},
  number={1},
  pages={133--172},
  year={2020}
}

@article{FHW94,
  title={5-chromatic {S}teiner triple systems},
  author={Fugere, Jean and Haddad, Lucien and Wehlau, David},
  journal={Journal of Combinatorial Designs},
  volume={2},
  number={5},
  pages={287--299},
  year={1994},
  publisher={Wiley Online Library}
}

@article{singer1938theorem,
  title={A theorem in finite projective geometry and some applications to number theory},
  author={Singer, James},
  journal={Transactions of the American Mathematical Society},
  volume={43},
  number={3},
  pages={377--385},
  year={1938},
  publisher={JSTOR}
}

@article{csima1984colouring,
  title={Colouring finite planes of odd order},
  author={Csima, J},
  journal={Discrete Mathematics},
  volume={49},
  number={3},
  pages={309},
  year={1984},
  publisher={Elsevier}
}

@article{knuth2000dancing,
  title={Dancing links},
  author={Knuth, Donald E},
  journal={arXiv preprint cs/0011047},
  year={2000}
}

@misc{prellberg2026constraintsatisfactionprogrammingnothreeinline,
      title={Constraint Satisfaction Programming for the No-three-in-line Problem}, 
      author={Thomas Prellberg},
      year={2026},
      eprint={2602.07751},
      archivePrefix={arXiv},
      primaryClass={math.CO},
      url={https://arxiv.org/abs/2602.07751}, 
}

@misc{v2026generalpositionproblemsurvey,
      title={The General Position Problem: A Survey (version 4)}, 
      author={Ullas Chandran S. V. and Sandi Klavžar and James Tuite},
      year={2026},
      eprint={2501.19385},
      archivePrefix={arXiv},
      primaryClass={math.CO},
      url={https://arxiv.org/abs/2501.19385}, 
}

@article{wood2004note,
  title={A note on colouring the plane grid},
  author={Wood, David R},
  journal={Geombinatorics},
  volume={13},
  number={4},
  pages={193--196},
  year={2004}
}

@article{hirschfeld2025arcs,
  title={Arcs, caps and generalisations in a finite projective space},
  author={Hirschfeld, James WP and Thas, Joseph A},
  journal={Mathematics},
  volume={13},
  number={9},
  pages={1489},
  year={2025},
  publisher={MDPI}
}

@article{Had99,
  title={On the chromatic numbers of {S}teiner triple systems},
  author={Haddad, Lucien},
  journal={Journal of Combinatorial Designs},
  volume={7},
  number={1},
  pages={1--10},
  year={1999},
  publisher={Wiley Online Library}
}

@book{arndt2010matters,
  title={Matters Computational: ideas, algorithms, source code},
  author={Arndt, J{\"o}rg},
  year={2010},
  publisher={Springer Science \& Business Media}
}

@book{hardy75,
  title={An introduction to the theory of numbers},
  author={Hardy, Godfrey Harold and Wright, Edward Maitland},
  year={1979},
  publisher={Oxford University Press}
}

@article{haralick1980increasing,
  title={Increasing tree search efficiency for constraint satisfaction problems},
  author={Haralick, Robert M and Elliott, Gordon L},
  journal={Artificial Intelligence},
  volume={14},
  number={3},
  pages={263--313},
  year={1980},
  publisher={Elsevier}
}

@article{flammenkamp1992progress,
  title={Progress in the no-three-in-line-problem},
  author={Flammenkamp, Achim},
  journal={Journal of Combinatorial Theory, Series A},
  volume={60},
  number={2},
  pages={305--311},
  year={1992},
  publisher={Elsevier}
}

@article{hall1975some,
  title={Some advances in the no-three-in-line problem},
  author={Hall, Richard R and Jackson, Terence H and Sudbery, Anthony and Wild, Ken},
  journal={Journal of Combinatorial Theory, Series A},
  volume={18},
  number={3},
  pages={336--341},
  year={1975},
  publisher={Elsevier}
}

@book{guy2004unsolved,
  title={Unsolved problems in number theory},
  author={Guy, Richard},
  volume={1},
  year={2004},
  publisher={Springer Science \& Business Media}
}

@book{kiss2019finite,
  title={Finite geometries},
  author={Kiss, Gyorgy and Sz\H{o}nyi, Tam{\'a}s},
  year={2019},
  publisher={Chapman and Hall/CRC}
}

@book{hirschfeld1998projective,
  title={Projective geometries over finite fields},
  author={Hirschfeld, James William Peter},
  year={1998},
  publisher={Oxford University Press}
}

@article{ebert1985partitioning,
  title={Partitioning projective geometries into caps},
  author={Ebert, Gary L},
  journal={Canadian Journal of Mathematics},
  volume={37},
  number={6},
  pages={1163--1175},
  year={1985},
  publisher={Cambridge University Press}
}

@article{lavrauw2014geometry,
  title={Geometry of the inversion in a finite field and partitions of $\text{PG}(2k- 1, q)$ in normal rational curves},
  author={Lavrauw, Michel and Zanella, Corrado},
  journal={Journal of Geometry},
  volume={105},
  number={1},
  pages={103--110},
  year={2014},
  publisher={Springer}
}

@book{godsil2001algebraic,
  title={Algebraic graph theory},
  author={Godsil, Chris and Royle, Gordon F},
  volume={207},
  year={2001},
  publisher={Springer Science \& Business Media}
}

@article{hall1974difference,
  title={Difference sets},
  author={Hall, M.},
  journal={Math. Centre Tracts.},
  volume={51},
  pages={1--26},
  year={1974}
}

@article{faina2002cyclic,
  title={The Cyclic Model for $\text{PG}(n,q)$ and a Construction of Arcs},
  author={Faina, Giorgio and Kiss, Gy{\"o}rgy and Marcugini, Stefano and Pambianco, Fernanda},
  journal={European Journal of Combinatorics},
  volume={23},
  number={1},
  pages={31--35},
  year={2002},
  publisher={Elsevier}
}

@book{scheinerman2011fractional,
  title={Fractional graph theory: a rational approach to the theory of graphs},
  author={Scheinerman, Edward R and Ullman, Daniel H},
  year={2011},
  publisher={Courier Corporation}
}

@article{kadelburg2005inequalities,
  title={Inequalities of {K}aramata, {S}chur and {M}uirhead, and some applications},
  author={Kadelburg, Zoran and Dukic, Dusan and Lukic, Milivoje and Matic, Ivan},
  journal={The Teaching of Mathematics},
  volume={8},
  number={1},
  pages={31--45},
  year={2005},
  publisher={Dru{\v{s}}tvo matemati{\v{c}}ara Srbije}
}

@book{qvist1952some,
  title={Some remarks concerning curves of the second degree in a finite plane},
  author={Qvist, Bertil},
  year={1952},
  publisher={Suomalainen tiedeakatemia}
}

@article{araujo2003daisy,
  title={Daisy structure in Desarguesian projective planes},
  author={Araujo-Pardo, Gabriela},
  journal={Journal of the Australian Mathematical Society},
  volume={74},
  number={2},
  pages={145--154},
  year={2003},
  publisher={Cambridge University Press}
}

@article{bacso2013twobloc,
  title={The 2-Blocking Number and the Upper Chromatic Number of $\text{PG} (2, q)$},
  author={Bacs{\'o}, G{\'a}bor and H{\'e}ger, Tam{\'a}s and Sz{\H{o}}nyi, Tam{\'a}s},
  journal={Journal of Combinatorial Designs},
  volume={21},
  number={12},
  pages={585--602},
  year={2013},
  publisher={Wiley Online Library}
}

@article{bacso2007upper,
  title={Upper chromatic number of finite projective planes},
  author={Bacs{\'o}, G{\'a}bor and Tuza, Zsolt},
  journal={Journal of Combinatorial Designs},
  volume={16},
  number={3},
  pages={221--230},
  year={2008},
  publisher={Wiley Online Library}
}

@article{akm2015balanupper,
  title={On the balanced upper chromatic number of cyclic projective planes and projective spaces},
  author={Araujo-Pardo, Gabriela and Kiss, Gy{\"o}rgy and Montejano, Amanda},
  journal={Discrete Mathematics},
  volume={338},
  number={12},
  pages={2562--2571},
  year={2015},
  publisher={Elsevier}
}

@article{bbmns21balaupper,
  title={On the balanced upper chromatic number of finite projective plane},
  author={Bl{\'a}zsik, Zolt{\'a}n L and Blokhuis, Aart and Miklavi{\v{c}}, {\v{S}}tefko and Nagy, Zolt{\'a}n L{\'o}r{\'a}nt and Sz{\H{o}}nyi, Tam{\'a}s},
  journal={Discrete Mathematics},
  volume={344},
  number={3},
  pages={112266},
  year={2021}
}

@article{araujo2019indaffine,
  title={On chromatic indices of finite affine spaces},
  author={Araujo-Pardo, Gabriela and Kiss and Gy{\"o}rgy and Rubio-Montiel Christian and V\'azquez \'Avila Adri\'an},
  journal={Ars Mathematica Contemporanea},
  volume={16},
  pages={67-79},
  year={2019}
}

@article{araujo2019inproy,
  title={On line colorings of finite projective spaces},
  author={Araujo-Pardo, Gabriela and Kiss, Gy{\"o}rgy and Rubio-Montiel, Christian and V{\'a}zquez-{\'A}vila, Adri{\'a}n},
  journal={Graphs and Combinatorics},
  volume={37},
  number={3},
  pages={891--905},
  year={2021},
  publisher={Springer}
}

\end{document}